\documentclass[12pt,french]{amsart}

\usepackage{amssymb}
\usepackage{graphicx}
\usepackage{epsf}
\usepackage{amsmath}
\usepackage[latin1]{inputenc}
\usepackage{amsfonts}
\usepackage{color}

\newtheorem{theo}{Theorem}
\newtheorem{prop}{Proposition}

\newtheorem{lemm}[theo]{Lemma}
\newtheorem{coro}[theo]{Corollary}
\newtheorem{defi}[theo]{Definition}

\newtheorem{rema}[theo]{Remark}

\def\cf{{\it cf. }}
\def\ie{{\it i.e. }}

\def\ti{$\sim$}

\def\CC{{\mathcal C}}

\def\shuffle{{\,\raise
1pt\hbox{$\scriptscriptstyle\cup{\mskip-4mu}\cup$}\,}}

\def\Young#1{\vbox{\smallskip\offinterlineskip
    \halign{&\vbox{##}\kern-\Thickness\cr #1}}}

\newdimen\Squaresize \Squaresize=12pt
\newdimen\Thickness \Thickness=.1pt
\newdimen\Correction \Correction=7pt

\def\Vide#1{\hbox{
       \vbox to \Squaresize{\vss
          \hbox to \Squaresize{\hss#1 \hss}\vss}
    \hskip-\Correction}
   \kern-\Thickness}

\def\Carre#1{\hbox{\vrule width \Thickness
   \vbox to \Squaresize{\hrule height \Thickness\vss
      \hbox to \Squaresize{\hss$\scriptstyle#1$\hss}
   \vss\hrule height\Thickness}
   \unskip\vrule width \Thickness}
   \kern-\Thickness}

\title{Keys and alternating sign matrices}

\author{Jean-Christophe Aval}\address[Jean-Christophe Aval]{LaBRI\\ Universit\'e Bordeaux 1\\ 351 cours
 de la Lib\'eration\\ 33405 Talence cedex\\ FRANCE}
\email{aval@labri.fr}
\urladdr{http://www.labri.fr/\ti aval/}

\date{}

\thanks{This work has been supported by the ANR (project MARS/06-BLAN-0193)} 

\begin{document} 

\maketitle 

\begin{abstract}
In \cite{key}, Lascoux and Sch\"utzenberger introduced a notion of  key associated to any Young tableau. More recently Lascoux defined the  key of an alternating sign matrix by recursively removing all $-1$'s in such matrices. But alternating sign matrices are in bijection with monotone triangles, which form a subclass of Young tableaux. We show that in this case these two notions of keys coincide. Moreover we obtain an elegant and direct way to compute the key of any Young tableau, and discuss consequences of our result.
\end{abstract}

\bigskip
\bigskip

\section{Introduction}
A  key is by definition a Young tableau whose columns are comparable for the inclusion order. They were introduced by A. Lascoux and M.-P. Sch\"utzenberger \cite{key,aretes} to study Demazure characters. For type $A$, irreducible characters (\ie Schur functions) are associated to all Young tableaux of a given shape. Demazure characters correspond to subsets of tableaux which can be described using keys (\cf \cite{lak, lit,mamede,mason} for recent works). 

Alternating sign matrices (ASM in short) are combinatorial objects that were extensively studied in the last two decades \cite{bressoud,propp}, with (at least) the great achievment of their enumeration \cite{zeil,kup}.
They may be seen as a representation of  square-ice configurations \cite{chern}, and surprisingly their numbers appeared in the context of the dense $O(n=1)$ or Temperley-Lieb loop model on the square grid \cite{RS1,RS2,saclay,gier}.

In \cite{chern}, A. Lascoux defines an operation which consists in iteratively  removing the $-1$'s of a given ASM to obtain a permutation matrix, called the key of the ASM. 
But ASM's are in bijection with a certain class of Young tableaux, called  monotone triangles (or Gog triangles in \cite{zeil}). Since a permutation (matrix) may be seen as a monotone triangle which is a key as a Young tableau, we may ask whether these two notions coincide. The answer, as stated without proof in \cite{chern}, is affirmative (\cf Corollary \ref{keylascoux}). 

More generally we extend Lascoux's operation of removing the $-1$'s to matrices associated to (unrestricted) Young tableaux and obtain (Theorem \ref{main}) a very simple way to compute the left key of a Young tableau, by far easier and quicker than the original definition \cite{key}, and even than the recent method presented in \cite{mason}.

This article is organized as follows. In Section 2, we recall the definitions of the left key of a Young tableau, introduce and prove our new way to compute it. We also explain how this method can be used for the computation of the right key. In Section 3, we examine consequences of our result for ASMs, including simple formulas for the number of ASMs with exactly one or two $-1$ entries.

\section{Keys of Young tableaux}\label{sec2}

In this paper, the French notation for tableau is used. A tableau may be seen as a product of columns: $T=C_1\,C_2\,\cdots\,C_l$, where a column is a strictly decreasing word (often identified with the set of its entries). For example, the tableau $T=
\Young{
\Carre{4}&&\cr
\Carre{2}&\Carre{5}&\cr
\Carre{1}&\Carre{2}&\Carre{5}\cr
}$ is the product of the 3 columns $C_1\, C_2\, C_3$ with $C_1=4\,2\,1$, $C_2=5\,2$ and $C_3=5$. Moreover, this tableau $T$ is a {\em Young tableau} because the entries are non-decreasing along its rows. The {\em shape} of a tableau is then the list $(H_1,\dots,H_l)$ of the heights of its columns. The shape of our tableau $T$ given above is $(3,2,1)$.

The {\em word} of a tableau is simply the result of its reading (column by column). The tableau $T$ reads: $421\ 52\ 5$ (the gaps are here to mark the changes of columns, and are of course unnecessary). We will now not distinguish a tableau and its word.

A {\em key} is a tableau such that for any $k\in\{1,\dots,l-1\}$, the column $C_{k+1}$ is a subword of the column $C_k$.

We recall that Schensted's insertion gives a bijection from words to pairs of Young tableaux $(P,Q)$ of same shape, with the second one standard. 
Two words are equivalent in the sense of the plactic congruence, characterized by Knuth, if they give the same tableau $P$ by insertion.
The following lemma (\cf \cite{aretes}) is the crucial step to associate a key to any Young tableau.

\begin{lemm}
Let $T$ be a Young tableau of shape $H=(H_1,\dots,H_l)$. Then for any permutation $I=(I_1,\dots,I_l)$ of $H$, there exists in the congruence class of  $T$ exactly one word $V=V_1\cdots V_l$, which is a product of columns of respective degrees $I_1,\dots,I_l$. If $J$ is another permutation of $I$ and $W=W_1\cdots W_l$ is the associated word, then $J_k\ge I_k$ implies that $V_k$ is a subword of $W_k$.
\end{lemm}

The words $V$ introduced in this lemma are called the {\em frank words} of $T$.
The keys are characterized by the fact that their columns (as words) commute (in the plactic monoid). 
For a tableau $T$, we shall say that $V_1$ is the {\em left factor} of the couple $(I,T)$, and $V_l$ the {\em right factor}. Thus, for each height $H_k$, we obtain this way a unique left (resp. right) factor of degree $H_k$. These columns are all ordered by inclusion. The {\em left (resp.right) key of the Young tableau} $T$ is then defined as the product of its left (resp. right) factors of respective sizes $H_1,H_2,\dots,H_l$.

In our example, the tableau $T$ is of shape $(3,2,1)$, thus there are six permutations of its shape, which correspond to the six frank words:

$$
\Young{
\Carre{4}&&\cr
\Carre{2}&\Carre{5}&\cr
\Carre{1}&\Carre{2}&\Carre{5}\cr
} \equiv
\Young{
\Carre{4}&\Carre{5}&\cr
\Carre{2}&\Carre{2}&\cr
&\Carre{1}&\Carre{5}\cr
} \equiv
\Young{
\Carre{4}&&\cr
\Carre{2}&&\cr
\Carre{1}&\Carre{5}&\Carre{5}\cr
&&\Carre{2}\cr
} \equiv
\Young{
\Carre{4}&\Carre{5}&\cr
&\Carre{2}&\Carre{5}\cr
&\Carre{1}&\Carre{2}\cr
} \equiv
\Young{
\Carre{4}&\Carre{5}&\Carre{5}\cr
&\Carre{2}&\Carre{2}\cr
&&\Carre{1}\cr
} \equiv
\Young{
\Carre{4}&&\cr
\Carre{2}&\Carre{5}&\Carre{5}\cr
&&\Carre{2}\cr
&&\Carre{1}\cr
} 
$$

We chose a planar presentation of words (\ie in form of tableaux) to insist on the fact that the frank words may be computed using the ``jeu de taquin''. We will use the ``jeu de taquin'' and refer to classical references \cite{taq} for definitions, and to \cite{fulton}, Appendix A.5 for its use in the computation of keys.
The left and right keys of our tableau $T$ are respectively  
$$\Young{
\Carre{4}&&\cr
\Carre{2}&\Carre{4}&\cr
\Carre{1}&\Carre{2}&\Carre{4}\cr
}\ \  {\rm and}\ \ 
\Young{
\Carre{5}&&\cr
\Carre{2}&\Carre{5}&\cr
\Carre{1}&\Carre{2}&\Carre{5}\cr
}
\, .$$

We observe that to compute the left key of a tableau $T$, it is sufficient to compute its left factors; thus we may restrict the computation to frank words relative to permutations $(I_k,I_1,I_2,\dots,I_{k-1},I_{k+1},\dots,I_l)$.

Now we associate to any tableau $T$ a sign matrix $M(T)$.
\begin{defi}\label{SM}
A {\em sign matrix} is a matrix $M=(M_{i,j})$ such that:
\begin{itemize}
\item $\forall i,j,\ M_{i,j}\in\{-1,0,1\}$;
\item $\forall i,j, \sum_{r=1}^i M_{r,j}\in\{0,1\}$;
\item $\forall i,j, \sum_{s=1}^j M_{i,s}\ge 0$.
\end{itemize}
\end{defi}

The bijection between Young tableaux and sign matrices is a generalization of the well-known bijection between ASM and monotone triangle (\cf \cite{zeil}).
We observe the apparitions and disparitions of the entries in the columns (from right to left), and we translate it matricially:
$$M(T)_{ij}=\left\{
\begin{array}{rl}
1& {\rm if\ }j \in C_{l-i+1} {\rm \ and}\ j\not\in C_{l-i+2}\cr
-1& {\rm if\ }j \not\in C_{l-i+1} {\rm \ and}\ j\in C_{l-i+2}\cr
0& {\rm if\ not}
\end{array}
\right.$$
with the convention that $C_{l+1}$ is empty.

\begin{prop}
The application that sends a Young tableau $T$ to a matrix $M(T)$ is a bijection from the set of Young tableaux with $m$ columns and entries in $\{1,\dots,n\}$ to sign matrices with $m$ rows and $n$ columns.
\end{prop}
\proof 
First, if we start with a young tableau $T$, the matrix obtained is a sign matrix because: 
\begin{itemize}
\item since $T$ is a Young tableau, the elements of a column of $T$ are all distinct, thus in a given column of $M(T)$, the non-zero entries start with a $1$, and then alternate, whence the second condition of sign matrices;
\item the elements of $T$ are weakly increasing along the rows, which translates matricially as the third condition of the sign matrices definition.
\end{itemize}

Conversely, if we start from a sign matrix, we construct a tableau, which is a Young tableau for the same reason.
\endproof

For example, the tableau $T=
\Young{
\Carre{5}&\Carre{5}&&\cr
\Carre{2}&\Carre{4}&\Carre{5}&\cr
\Carre{1}&\Carre{2}&\Carre{4}&\Carre{6}\cr
}$
is associated to the sign matrix 
$$M(T)=\left[
\begin{array}{cccccc}
0&0&0&0&0&1\cr
0&0&0&1&1&-1\cr
0&1&0&0&0&0\cr
1&0&0&-1&0&0\cr
\end{array}
\right]$$

A simple observation is that a Young tableau is a key if and only if its sign matrix does not contain any $-1$.
Since our goal is to associate to any Young tableau its (left) key, we introduce a way to remove the $-1$'s.
This is done through an {\em elimination process}. This process, as we shall see in the next section, is an extension of the {\em removing process} defined in \cite{chern} for monotone triangles, to general Young tableaux. It should be observed here that this process is more than just the restriction of Lascoux's process to the sub-quadrant of an ASM, as may be seen on the previous example, where the first two rows are clearly not a part of an ASM.

A $-1$ entry in a matrix $M$ is said {\em removable} if there is no $-1$ in the rows above it, nor in its row and on its right, \ie for the $-1$ in position $(a,b)$ (we use matricial coordinates) if:
$$\forall i<a,\ \forall j,\ M_{i,j}\neq -1,\ \ {\rm and}\ \forall j>b,\ M_{a,j}\neq-1.$$
For such a $-1$ in a given matrix, its {\em neighbours} are the entries $M_{i,j}$ equal to $1$ such that
$$i\le a,\ j\le b,\ {\rm and}\ \forall i\le k\le a,\ \forall j\le l\le b,\ M_{k,l}\neq 1,$$
\ie the rectangle of South-East corner $(a,b)$ and of North-West corner $(i,j)$ contains no other entry equal to $1$ than the neighbour itself. For a given removable $-1$, the union of these rectangles is a Ferrers diagram:
$$\begin{array}{cccccccc}
.&.&.&.&.&.&.&1\cr
.&.&.&.&.&1&0&0\cr
.&.&.&.&.&0&0&0\cr
.&.&.&1&0&0&0&0\cr
.&.&.&0&0&0&0&0\cr
.&.&.&0&0&0&0&0\cr
1&\theta&\theta&\theta&\theta&\theta&\theta&-1\cr
\end{array}$$
where the $\theta$ entries are either $0$'s or $-1$'s.
To {\em remove} a (removable) $-1$ consists in replacing it by a $0$, and to replace its $n$ neighbours by $0$ and to place $n-1$ entries equal to $1$ such that they form a new Ferrers diagram whose inner corners are the former $n$ neighbours. To make this definition unambiguous, we precise that in our context the {\em inner corners} of a Ferrers diagram are the dots in the following picture.

\centerline{\epsffile{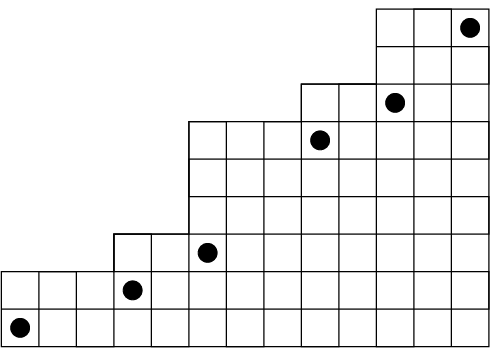}}
\smallskip

 On the example, this means replace the sub-matrix given above by:
$$\begin{array}{cccccccc}
.&.&.&.&.&{\bf 1}&.&{\bf 0}\cr
.&.&.&{\bf 1}&.&{\bf 0}&0&0\cr
.&.&.&.&.&0&0&0\cr
{\bf 1}&.&.&{\bf 0}&0&0&0&0\cr
.&.&.&0&0&0&0&0\cr
.&.&.&0&0&0&0&0\cr
{\bf 0}&\theta&\theta&\theta&\theta&\theta&\theta&{\bf 0}\cr
\end{array}$$

Now our main result is the following;

\begin{theo}\label{main}
Let $T$ be a Young tableau, and $M(T)$ the sign matrix associated to $T$. By removing all the $-1$'s of $M(T)$ by the process described, we obtain a sign matrix associated to a tableau $U$ which is the left key of $T$.
\end{theo}

Let us deal with our example $T=
\Young{
\Carre{5}&\Carre{5}&&&\cr
\Carre{2}&\Carre{4}&\Carre{5}&\Carre{6}&\cr
\Carre{1}&\Carre{2}&\Carre{4}&\Carre{4}&\Carre{6}\cr
}$. Its sign matrix is \label{page5}
$$M(T)=\left[
\begin{array}{cccccc}
0&0&0&0&0&1\cr
0&0&0&1&0&0\cr
0&0&0&0&1&-1\cr
0&1&0&0&0&0\cr
1&0&0&-1&0&0\cr
\end{array}
\right]$$

Now we compute:
$$M(T)\ \longrightarrow\ \left[
\begin{array}{cccccc}
0&0&0&0&1&0\cr
0&0&0&1&0&0\cr
0&0&0&0&0&0\cr
0&1&0&0&0&0\cr
1&0&0&-1&0&0\cr
\end{array}
\right]
\ \longrightarrow\ 
\left[
\begin{array}{cccccc}
0&0&0&0&1&0\cr
0&1&0&0&0&0\cr
0&0&0&0&0&0\cr
1&0&0&0&0&0\cr
0&0&0&0&0&0\cr
\end{array}
\right]=M(U)$$
with $U=\Young{
\Carre{5}&\Carre{5}&&&\cr
\Carre{2}&\Carre{2}&\Carre{5}&\Carre{5}&\cr
\Carre{1}&\Carre{1}&\Carre{2}&\Carre{2}&\Carre{5}\cr
}$, the left key of $T$.

\medskip

The proof of the theorem starts with the following lemma.

\begin{lemm}\label{lemN1}
When we apply the ``jeu de taquin'' to two-column tableaux to compute frank words:

\centerline{\epsffile{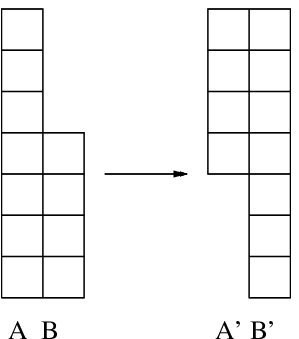}}

\noindent it is clear that if $b\in B\cap A$ then $b\in B'\cap A'.$ Now if $A$ and $B$ do not commute, \ie if $B-A\neq\emptyset$, let $\tilde b=max(B-A)$ and $\tilde a=max(\{a\in A-B,\ a\le \tilde b\}).$ We have
\begin{enumerate}
\item $\tilde a\in A'$;
\item let $\tilde B$ be the column obtained by replacing $\tilde b$ by $\tilde a$ in $B$ (maybe reordered), then if $A\tilde B\longrightarrow \tilde A'\tilde B'$, we have: $\tilde A'=A'$.
\end{enumerate}
\end{lemm}

\proof
For the first point, if we suppose  that $\tilde a$ slips to $B'$ during the ``jeu de taquin'', $\tilde a$ should have the hole on its right (and obviously below $\tilde b$). Let $\{x_1,x_2,\dots,x_p\}$ be the entries in the left column above $\tilde a$ and not above $\tilde b$:
$$\Young{
\Carre{x_1}&\Carre{\tilde b}\cr
\Carre{x_2}&\Carre{}\cr
\Carre{\cdots}&\Carre{}\cr
\Carre{x_p}&\Carre{}\cr
\Carre{\tilde a}&\Carre{X}\cr
}$$
where the $X$ is the hole. Since $\forall i\in\{1,\dots,p\},\ \tilde a<x_i\le \tilde b$ and by definition of $\tilde a$, we have that all the $x_1,\dots,x_p$ are in the right column. But these entries should be below the box that contains $\tilde b$ (since $x_i\le \tilde b$) and above the hole (since $x_i\ge \tilde a$). This is impossible since we only have here $p-1$ boxes.

For the second point, we observe that for any element $a\in A$, the fact that $a$ stays in $A'$ (or not) is the same whether $B$ contains $\tilde b$ or $\tilde a$:
\begin{itemize}
\item for $a>\tilde b$, nothing changes;
\item the case $a=\tilde b$ does not occur;
\item if $\tilde a <a< \tilde b$, $a$ is in $A$ if and only if $a$ is in $B$ (or $\tilde B$), thus this is clear;
\item for $a=\tilde a$, this is the first point of the lemma;
\item for $a<\tilde a$, nothing changes.
\end{itemize}
\endproof

\noindent{\it Proof of the Theorem \ref{main}.}

Since a sign matrix without any $-1$ is the matrix of a key, we only have to show that the elimination of a $-1$ does not change the left key.

Thus let $T$ be a tableau and $T'$ the tableau resulting from the elimination of a $-1$ (its removable $-1$).
This elimination concerns an entry $a$ appeared in column $j$ and disappeared in column $i$ (we recall that this story reads from right to left).

The elimination of the $-1$ consists in replacing in $T$ the entry $a$ in the column $k$ (for $i+1\le k\le j$) by the greatest of all entries of $C_i$ less or equal to $a$ and that are not in $C_k$. 

Now let us compare the computation of frank words for $T$ and $T'$. For the left factors of size $\ge H_{i}$, nothing has changed from  $T$ to $T'$.
Let $k\ge i+1$. By definition of the removable $-1$, there is no $-1$ in the rows of $M(T)$ above the row $i$. Thus the columns on the right of $C_i$ commute, whether we consider $T$ or $T'$. We have to deal with the action of the ``jeu de taquin'' on $C_iC_k$. In the case of $T$, we know from Lemma \ref{lemN1} that the entry $a$ is replaced by the greatest of all entries of $C_i$ less or equal to $a$ and which are not in $C_k$. Moreover, for all entries less than $a$, the computation in the same for $T$ and $T'$.

Thus $T$ and $T'$ have the same left key.
\qed

\smallskip

Now we show how this process can be used to compute the right key of a Young tableau. As we shall see in Theorem \ref{teocomp}, we can reduce it to the computation of a left key.

Let us introduce the notion of complement of a tableau, whose definition was already given in \cite{mono}.

\begin{defi}
Let $T$ be a tableau with entries in $\{1,2,\dots,n\}$ and columns $T=C_1,\dots,C_l$. The {\em complement} of $T$, denoted by $\CC(T)$ is defined as $\CC(T)=D_1,\dots,D_l$ with for $1\le j\le n$ and $1\le i\le l$:
$$j\in D_i\ \Longleftrightarrow \ j\not\in C_{l-i+1}.$$
\end{defi}

For example if 
$T=\Young{
\Carre{5}&\Carre{5}\cr
\Carre{2}&\Carre{4}&\Carre{5}&\Carre{6}\cr
\Carre{1}&\Carre{2}&\Carre{4}&\Carre{4}&\Carre{6}\cr}
$ then $\CC(T)=
\Young{
\Carre{5}\cr
\Carre{4}&\Carre{5}&\Carre{6}\cr
\Carre{3}&\Carre{3}&\Carre{3}&\Carre{6}&\Carre{6}\cr
\Carre{2}&\Carre{2}&\Carre{2}&\Carre{3}&\Carre{4}\cr
\Carre{1}&\Carre{1}&\Carre{1}&\Carre{1}&\Carre{3}\cr}.$

\begin{rema}
The complement of a tableau $T$ depends of course on the set of entries that we consider. For example, if we had considered that our previous example $T$ has entries in $\{1,\dots,7\}$, the only difference would have been the presence of a box with a $7$ in each column of $\CC(T)$. Thus it should be clear that the next results (Proposition \ref{compyoung}, Lemma \ref{lemstar} and Theorem \ref{teocomp}) remain true if we add some entry to the set of entries.
\end{rema}

\begin{prop}\label{compyoung}
If $T$ is a Young tableau, then $\CC(T)$ is also a Young tableau.
\end{prop}

\proof
This is a direct consequence of the following Lemma \ref{lemstar}.
\endproof

\begin{lemm}\label{lemstar}
Let $T$ be a tableau with $l$ columns and $M$ its sign matrix. If $M'$ is the sign matrix of its complement, then we have:
$$\forall i>1,\ \forall j,\ M_{i,j}=M'_{l-i+1,j}.$$
\end{lemm}
\proof
Let us introduce the columns: $T=C_1\dots C_l$ and $\CC(T)=D_1\dots D_l$. For $i>1$, we may write:
\begin{eqnarray*}
M_{i,j}=1 & \Longleftrightarrow & j\in C_{l-i+1}-C_{l-i+2}\\ 
 & \Longleftrightarrow & j\in D_i-D_{i-1}\\
 & \Longleftrightarrow & M'_{l-i+1,j}=1\\
\end{eqnarray*}
and the same holds for $-1$'s, which proves the lemma.
\endproof

\begin{theo}\label{teocomp}
Let $T$ be a Young tableau. The complement of the right key of $T$ is the left key of the complement of $T$.
\end{theo}
\proof
A first observation is that we can restrict to Young tableaux with exactly two columns, since frank words are computed by iteratively applying the jeu de taquin to pairs of columns.

Now a consequence of Lemma \ref{lemN1} is that the result of the action of the ``jeu de taquin'' $AB\longrightarrow A'B'$ depends only on $A\cap B=A'\cap B'$, $A-B$ and $B-A$, which we may write (where the symbol $\sqcup$ means ``disjoint union''):
$$A'=(A\cap B) \sqcup T_1(A-B,B-A)$$
$$B'=(A\cap B) \sqcup T_2(A-B,B-A)$$
with $T_1(A-B,B-A)\cap T_2(A-B,B-A)=\emptyset$.

If we denote $CD=\CC(AB)$, we have:
$$C-D=A-B\ \ {\rm and}\ \ D-C=B-A.$$

We may write
$$C'=(C\cap D) \sqcup T_1(C-D,D-C)$$
$$D'=(C\cap D) \sqcup T_2(C-D,D-C)$$
and we want to check that $A'\sqcup D'=\{1,\dots,n\}=B'\sqcup C'$. By symmetry, we restrict to the first equality which is proved as follows:
$$\left\{
\begin{array}{rcl}
A'\cap D'&=&\Big((A\cap B)\sqcup T_1(A-B,B-A\Big)\cap\Big((C\cap D) \sqcup T_2(C-D,D-C)\Big)\cr
&=&\Big((A\cap B)\sqcup T_1(A-B,B-A\Big)\cr
&&$\ \ \ \ \ \ \ \ \ \ \ \ \ \ \ \ \ \ \ \ \ \ $\cap\Big(\big(\{1,\dots,n\}-(A\cup B)\big) \sqcup T_2(A-B,B-A)\Big)\cr
&=&\emptyset\cr
|A'|+|D'|&=&n.
\end{array}
\right.$$
\endproof

We deduce from Theorem \ref{teocomp} an easy way to compute the right key of a Young tableau by taking the complement, applying Theorem \ref{main} to obtain the left key, then taking the complement again. For example, if we consider $T=\Young{
\Carre{5}&\Carre{5}\cr
\Carre{2}&\Carre{4}&\Carre{5}&\Carre{6}\cr
\Carre{1}&\Carre{2}&\Carre{4}&\Carre{4}&\Carre{6}\cr}
$, we obtain easily that the left key of its complement $\CC(T)=
\Young{
\Carre{5}\cr
\Carre{4}&\Carre{5}&\Carre{6}\cr
\Carre{3}&\Carre{3}&\Carre{3}&\Carre{6}&\Carre{6}\cr
\Carre{2}&\Carre{2}&\Carre{2}&\Carre{3}&\Carre{4}\cr
\Carre{1}&\Carre{1}&\Carre{1}&\Carre{1}&\Carre{3}\cr}$
is $\Young{
\Carre{5}\cr
\Carre{4}&\Carre{5}&\Carre{5}\cr
\Carre{3}&\Carre{3}&\Carre{3}&\Carre{5}&\Carre{5}\cr
\Carre{2}&\Carre{2}&\Carre{2}&\Carre{3}&\Carre{3}\cr
\Carre{1}&\Carre{1}&\Carre{1}&\Carre{1}&\Carre{1}\cr}$
thus its right key is $\Young{
\Carre{6}&\Carre{6}\cr
\Carre{4}&\Carre{4}&\Carre{6}&\Carre{6}\cr
\Carre{2}&\Carre{2}&\Carre{4}&\Carre{4}&\Carre{6}\cr}
$ .

\section{Keys of alternating sign matrices}

Let us switch to the origin of this work, through A. Lascoux's paper \cite{chern}, that is to the context of alternating sign matrices. In this paper, using the notion of key, A. Lascoux obtains a description of Grothendieck polynomials (he already obtained a description of Schubert polynomials in \cite{schub}). In most questions where ASM's appear in algebraic combinatorics or theoretical physics, weights are assigned to ASM's. One of these weights can be specialized to $1$, which leads to enumeration, and the end of this third section deals with enumeration. But this is not the case for the weight giving Grothendieck polynomials.

An {\em alternating sign matrix} (ASM in short) of size $n$ is a square matrix with entries in $\{-1,0,1\}$ such that along each row and column:
\begin{itemize}
\item the sum of the entries is equal to $1$;
\item the non-zero entries alternate in sign.
\end{itemize}

It is clear that an ASM is in particular a sign matrix (in the sense defined in the previous section). Let us give an example of an ASM of size 5:
$$\left(
\begin{array}{ccccc}
0&1&0&0&0\cr
0&0&1&0&0\cr
1&-1&0&1&0\cr
0&1&-1&0&1\cr
0&0&1&0&0
\end{array}
\right)$$

The bijection defined in Section \ref{sec2} between sign marices and Young tableaux specializes here in the well-known (\cite{chern}) bijection between ASM's and monotone triangles. A {\em monotone triangle} is by definition a Young tableau of shape $(n,n-1,\dots,2,1)$ with entries in $\{1,2,\dots,n\}$ and such that the entries are non-decreasing along each diagonal (from South-East to North-West in French notation). As an example, we give the monotone triangle associated to the ASM given above:
$$\Young{
\Carre{5}\cr
\Carre{4}&\Carre{5}\cr
\Carre{3}&\Carre{4}&\Carre{4}\cr
\Carre{2}&\Carre{2}&\Carre{3}&\Carre{3}\cr
\Carre{1}&\Carre{1}&\Carre{1}&\Carre{2}&\Carre{2}\cr}$$

In \cite{chern}, A. Lascoux defined a process of elimination of the $-1$'s of an ASM, to obtain an ASM without any $-1$, \ie the matrix of a permutation, which Lascoux call ``the key of the ASM''. We may temporarily use the notions of matrix-key (or M-key) to refer to the key defined matricially by Lascoux, and of tableau-key (or T-key) to refer to the key defined originally on tableaux by Lascoux and Sch\"utzenberger in \cite{key}. Lascoux states without proof in \cite{chern} that the bijection between ASM and monotone triangles exchanges these two notions. This assertion now becomes a corollary of Theorem \ref{main}.

\begin{coro}\label{keylascoux}
Let $M$ be an ASM and $T$ its associated monotone triangle. If $K$ and $U$ are repectively the M-key of $M$ and the T-key of $T$, then $U$ is the monotone triangle in bijection with $K$.
\end{coro}

\proof
This is a direct consequence of Theorem \ref{main}, since the M-key is precisely computed by matricial rules defined in Section \ref{sec2}.
\endproof

Thus we may simply speak about ``key'' for ASM's or monotone triangles. We mention here that to index Schubert cells,
Ehressmann \cite{eh} used combinatorial objects which can be identified to our keys,
and showed that the order on keys corresponds to the natural order on cells
(later called Bruhat order).

\medskip

Now for combinatorial reasons, and to easily obtain the number of ASM's with exactly one or two $-1$'s, we introduce the notion of pseudo-key of an ASM.
The definition is almost the same as the definition of the key, but the process of elimination of the $-1$'s is simpler.

Here we look at a $-1$ in position $(i,j)$ in a given ASM $M$ such that there is no other $-1$ in its North-West quadrant, \ie $$\forall k\le i, l\le j, (k,l)\neq(i,j), M_{k,l}\neq -1.$$

Now in the matrix $M$ we simply replace the pattern 

$$\begin{array}{ccc}
0&\cdots&1\cr
\vdots&&\vdots\cr
1&\cdots&-1
\end{array}$$ 
by

$$\begin{array}{ccc}
1&\cdots&0\cr
\vdots&&\vdots\cr
0&\cdots&0
\end{array}.$$ 

The {\em pseudo-key} $pK(M)$ of an ASM $M$ is obtained by iteratively removing all $-1$'s according to this process. It is quite easy to check that the result does not depend on the order in which we perform these eliminations. It is clear that the pseudo-key can be defined of any sign matrix, and in the case of the tableau given page \pageref{page5}, the matricial computation is the following, which shows that the pseudo-key is in general different from the key.

$$M(T)=\left[
\begin{array}{cccccc}
0&0&0&0&0&1\cr
0&0&0&1&0&0\cr
0&0&0&0&1&-1\cr
0&1&0&0&0&0\cr
1&0&0&-1&0&0\cr
\end{array}
\right]
\longrightarrow\ 
\left[
\begin{array}{cccccc}
0&0&0&0&1&0\cr
1&0&0&0&0&0\cr
0&0&0&0&0&0\cr
0&1&0&0&0&0\cr
0&0&0&0&0&0\cr
\end{array}
\right]=pK(M(T))$$

Let us say a word about the lattice structure on ASM's. This notion comes from bijection with monotone triangles. It is quite easy to define the supremum (resp. infimum) of a family of monotone triangles, by just taking the supremum of the numbers in each box composing the monotone triangles. The {\em lattice structure} on ASM's (which extends Bruhat's order) is then inhereted from the lattice structure on monotone triangles. We mention the following fact, whose proof is straightforward.

\begin{prop}
Let $M$ be an ASM, $K(M)$ its key and $pK(M)$ its pseudo-key. For the order in the lattice of ASM's, we have:
$$pK(M)\le K(M)\le M.$$
\end{prop}

Now we conclude our work with two enumerative results. Let us denote by $A_n^{(k)}$ the number of alternating sign matrices with exactly $k$ entries equal to $-1$ (for example, $A_n^{(0)}=n!$).

\begin{prop}
We have for $n\ge 3$:
$$A_n^{(1)}=\frac{(n!)^2}{(3!)^2(n-3)!}$$
and for $n\ge 6$:
$$A_n^{(2)}=(n!)^2\Big(\frac{1}{2592(n-6)!}+\frac{11}{3600(n-5)!}+\frac{1}{288(n-4)!}\Big).$$
\end{prop}
\proof
We will only show how to prove the first result, the second result is proved in the same manner.

The proof will be in two steps: first we will prove that $A_n^{(1)}$ is equal to the total number of patterns $132$ in all permutations of $S_n$, then we will prove that this number is precisely $\frac{(n!)^2}{(3!)^2(n-3)!}$.

To start, we consider an ASM $M$ with exactly one $-1$. We compute its pseudo-key $pK(M)$ which is a permutation matrix. This operation gives a pattern $132$ in the permutation $\sigma$ associated to the matrix $pK(M)$, \ie a triple $(i,j,k)$ such that:
$$i<j<k,\ \ {\rm and}\ \ \sigma(i)<\sigma(k)<\sigma(j).$$
This fact is illustrated as:
$$\begin{array}{ccccc}
&&1&&\cr
&&\vdots&&\cr
1&\cdots&-1&\cdots&1\cr
&&\vdots&&\cr
&&1&&\cr
\end{array}
\longrightarrow
\begin{array}{ccccc}
{\bf 1}&&0&&\cr
&&\vdots&&\cr
0&\cdots&0&\cdots&{\bf 1}\cr
&&\vdots&&\cr
&&{\bf 1}&&\cr
\end{array}
$$
where the three remaining $1$'s form the pattern $132$.
If we keep track of this triple $(i,j,k)$, we obtain a bijection between ASM's with exactly one $-1$ and permutations with a marked pattern $132$. This is the first point of our proof.

Next, to compute the total number of patterns $132$ in all permutations of $S_n$, we observe that a pattern $132$ is given by:
\begin{itemize}
\item the choice of a triple $1\le i<j<k\le n$, which gives $n\choose3$ possibilities;
\item the choice of a triple $1\le \sigma(i)<\sigma(k)<\sigma(j)\le n$, which gives $n\choose3$ possibilities;
\item the choice of any permutation in the remaining $n-3$ positions.
\end{itemize}
Thus the total number of patterns $132$ is ${n\choose3}^2\times(n-3)!=\frac{(n!)^2}{(3!)^2(n-3)!}$, which completes the proof.
\endproof



\end{document}